\documentclass{article}
 
  \newtheorem{thm}{Theorem}[section]
 
  \newtheorem{lem}[thm]{Lemma}
  
  \newtheorem{prop}[thm]{Proposition}
  
\newtheorem{prob}{Problem}
 
\begin{document}
 
  \title{Computable structures of rank $\omega_{1}^{CK}$}
 
  \author{J.\ F.\ Knight and J.\ Millar}
 \date{}
  \maketitle
 
  \section{Introduction}

We consider only computable languages and countable structures.  In
measuring complexity, we identify a structure with its atomic diagram. 
Makkai \cite{M}  produced an example of an arithmetical structure of
Scott rank $\omega_{1}^{CK}$.  Here we obtain a computable structure  of
Scott rank $\omega_{1}^{CK}$ in two different ways.  The first way is to
start with Makkai's example and code it in a computable structure in a way
that preserves the rank.  We give a general coding procedure that
transforms any hyperarithmetical structure $\mathcal{A}$ into a
computable structure $\mathcal{A}^*$ such that the rank of $\mathcal{A}$
is $\omega_1^{CK}$, or $\omega_1^{CK}+1$, or $<\omega_1^{CK}$ iff the same
is true of $\mathcal{A}^*$.  The same coding procedure has been used
recently in \cite{G-H-K-M-M-S} for several further results. 

The second way that we obtain a computable structure with Scott rank
$\omega_1^{CK}$ is to examine Makkai's construction, and note that with a
little more care, the resulting structure of Scott rank $\omega_1^{CK}$
can be made computable.  We recall a construction, described nicely in \cite{Mo}, 
which associates trees with structures in such a way that a computable tree 
gives rise to a computable structure, and the tree has a path just in case the
corresponding structure has a non-trivial automorphism.  We show that if 
a computable tree is ``thin'', and has a path, but no hyperarithmetical one, 
then the resulting structure will have Scott rank $\omega_1^{CK}$.  
Having done this, we use an idea of Sacks to obtain a tree with the desired properties.  

The coding method is described in Section 2.  In Sections 3 and 4, we
re-work Makkai's construction.  Section 3 gives the passage from trees to
structures, and Section 4 constructs computable thin trees.  In the
remainder of the present section, we give a little background on Scott
rank.  

\subsection{Scott rank}

Scott \cite{Sc} proved the following.

\begin{thm} [Scott Isomorphism Theorem]
\label{thm1.1}

If $\mathcal{A}$ is a countable structure for a countable language, there
is an $L_{\omega_1,\omega}$-sentence $\sigma$ whose countable models are
just the isomorphic copies of $\mathcal{A}$.

\end{thm}   

In the proof of Theorem \ref{thm1.1}, Scott assigned countable ordinals
to the tuples in a countable structure $\mathcal{A}$, and to the structure
$\mathcal{A}$ itself.  He started with the following family of equivalence
relations $\equiv^\beta$ on tuples in $\mathcal{A}$.

\begin{enumerate}  

\item  $\overline{a}\equiv^0\overline{b}$ if $\overline{a}$ and
$\overline{b}$ satisfy the same quantifier-free formulas,

\item  for $\beta > 0$, $\overline{a}\equiv^\beta\overline{b}$ if for all
$\gamma < \beta$, for all $c\in\mathcal{A}$, there exists
$d\in\mathcal{A}$ such that
$\overline{a},c\equiv^{\gamma}\overline{b},d$ and for all
$d\in\mathcal{A}$, there exists $c\in\mathcal{A}$ such that
$\overline{a},c\equiv^\gamma\overline{b},d$.     

\end{enumerate}
\textbf{Remark}:  We could split 2 into two parts.  

\bigskip

(a)  For $\beta = \gamma + 1$, 
$\overline{a}\equiv^\beta\overline{b}$ if for all $c$ there exists $d$,

\hspace{.2in} and for all $d$ 
there exists $c$, such that $\overline{a},c\equiv^\gamma\overline{b},d$.  

\bigskip

(b)  For limit $\beta$, $\overline{a}\equiv^\beta\overline{b}$ if for all $\gamma <\beta$, $\overline{a}\equiv^\gamma\overline{b}$.  

\bigskip

If $\mathcal{A}$ is a countable structure, then for each tuple
$\overline{a}$, there is some countable ordinal $\beta$ such that all
$\overline{b}$ satisfying $\overline{a}\equiv^\beta\overline{b}$ are in
the orbit of $\overline{a}$.  The \emph{Scott rank} of the tuple $\overline{a}$, 
denoted by $SR(\overline{a})$, is the least
such $\beta$.  The \emph{Scott rank} of the structure $\mathcal{A}$, 
denoted by $SR(\mathcal{A})$, is the least ordinal $\alpha$
such that for all tuples $\overline{a}$ in $\mathcal{A}$,
$SR(\overline{a}) < \alpha$.  

\bigskip

We are interested in computable 
structures and their ranks.  We shall make use of  \emph{computable 
infinitary formulas} to state criteria for a computable structure to have 
various ranks.  Roughly speaking, computable infinitary formulas 
are infinitary formulas in which the 
infinite disjunctions and conjunctions are over c.e.\ sets.  For more 
information about these formulas, and for proofs of various facts stated 
below, see \cite{A-K}.                 
If $\mathcal{A}$ is a computable (or hyperarithmetical) structure, 
then for any tuple $\overline{a}$ in
$\mathcal{A}$ and any computable ordinal
$\beta$, there is a computable infinitary formula defining the
$\equiv^\beta$-class of $\overline{a}$.  Using Barwise-Kreisel Compactness, we
can show that the family of finite partial functions preserving
$\equiv^{\omega_1^{CK}}$ has the back-and-forth property.  It follows
that if $\mathcal{A}$ is hyperarithmetical, then 
$SR(\mathcal{A})\leq\omega_1^{CK}+1$.   

\begin{prop}
\label{prop1.2}

Suppose $\mathcal{A}$ is a hyperarithmetical structure.  
\begin{enumerate}

\item  If for some tuple $\overline{a}$, the orbit is not
defined by any computable infinitary formula, then $SR(\mathcal{A}) =
\omega_1^{CK}+1$.

\item  If for some computable ordinal $\alpha$, the orbits of all
tuples are defined by computable $\Pi_\alpha$ formulas, then
$SR(\mathcal{A}) < \omega_1^{CK}$.  

\item  Otherwise, $SR(\mathcal{A}) = \omega_1^{CK}$. 
\end{enumerate}

\end{prop} 

There are familiar examples of computable structures of Scott
rank $\omega_1^{CK}+1$.  Harrison \cite{H} showed that there is
a computable ordering of type $\omega_1^{CK}(1+\eta)$.  In the Harrison
ordering, the orbit of an element outside the well-ordered initial
segment of type $\omega_1^{CK}$ is not defined by any computable
infinitary formula, so the ordering has Scott rank $\omega_1^{CK}+1$. 
There are related examples of Boolean algebras and Abelian $p$-groups which also
have Scott rank $\omega_1^{CK}+1$ \cite{G-H-K-S}.

Makkai \cite{M} showed the following.

\begin{thm}[Makkai]
\label{thm1.3}

There exists an arithmetical structure $\mathcal{A}$ such that\linebreak 
$SR(\mathcal{A}) = \omega_{1}^{CK}$.

\end{thm}

In the present paper, we give two different proofs of the following.  

\begin{thm}
\label{thm1.4}

There exists a computable structure $\mathcal{A}$ such that $SR(\mathcal{A}) =
\omega_1^{CK}$.

\end{thm}

\section{Coding}

The first proof of Theorem \ref{thm1.4} takes Theorem \ref{thm1.3} as
given, and combines it with the result below.  
 
  \begin{thm}
\label{thm2.1}
 
Suppose $\mathcal{A}$ is a hyperarithmetical structure.  Then there 
is a computable structure $\mathcal{A}^*$, with a copy of
$\mathcal{A}$ definable in $\mathcal{A}^*$ by computable infinitary
formulas, such that
 
  \begin{enumerate}
 
  \item  If $\mathcal{A}$ has computable rank, then so does $\mathcal{A}^*$.
 
  \item  If $\mathcal{A}$ has rank $\omega_{1}^{CK}$, then so does
  $\mathcal{A}^*$.
 
  \item  If $\mathcal{A}$ has rank $\omega_{1}^{CK}+1$, then so does
  $\mathcal{A}^*$.
 
  \end{enumerate}
 
  \end{thm}

In what follows, we first describe the structure $\mathcal{A}^*$.  Then, we review
some basic notions and results, before proving that the ranks of
$\mathcal{A}$ and $\mathcal{A}^*$ are related as in Theorem \ref{thm2.1}.   

\subsection{Constructing $\mathcal{A}^*$}  

We suppose that $\mathcal{A}$ is a hyperarithmetical structure.  We may suppose that the universe is computable and the language is relational.  For each relation symbol $P$ from the language of $\mathcal{A}$, and each tuple $\overline{a}$ from $\mathcal{A}$, having the arity of $P$, the pair $(P,\overline{a})$ corresponds to an infinite set
  of points in the structure $\mathcal{A}^*$.  On this set, we put a copy
  of one of a pair of structures $\mathcal{C}_{1}$,
  $\mathcal{C}_{2}$, depending on whether $\mathcal{A}\models
  P\overline{a}$.  The structures that we choose will be linear
orderings.  

We give a precise description of the
structure $\mathcal{A}^*$, in terms of these linear orderings.  There are
unary relations $L$, $L^*$, $A$, $U$, and $T$, and binary relations $E$ and
  $<$.  The unary relations are disjoint, and their union is the universe
  of $\mathcal{A}^*$.  The elements of $L$ represent symbols in the 
  language
  of $\mathcal{A}$.  If $r\in L$ represents the $i^{th}$ 
  relation symbol, then $r$ is connected via $E$ to a cycle of length
$i+1$ in 
  $L^*$.
  The cycles are disjoint and exhaust $L^*$.  The elements of $A$ 
  represent
  elements of $\mathcal{A}$.  

For each $n$-tuple $a_{0},\ldots,a_{n-1}$ from
  $\mathcal{A}$, there is a corresponding $n$-tuple $u_{0},\ldots,u_{n-1}$
  from $U$, where $E$ connects $u_{i}$ to $a_{i}$.  In
  addition $E$ takes $u_{i}$ to $u_{i+1}$, and $u_{n-1}$ to $u_{0}$, 
  forming a cycle.  The cycles in $U$ are all distinct, and they exhaust
$U$.
  For each $r\in L$, representing an $n$-ary relation symbol $R$, and each
  $n$-tuple $\overline{a} = (a_{0},\ldots,a_{n-1})$ in $\mathcal{A}$, if,
  $u_{0}$ is the element connected to $a_{0}$ in in the corresponding 
  cycle
  in $U$, then $u_{0}$ and $r$ are both connected via $E$ to each element 
  in
  an infinite subset $T_{(R,\overline{a})}$ of $T$.  The sets
  $T_{(R,\overline{a})}$ are disjoint and they exhaust $T$.  There is one
more binary relation on $T$, which is the union of the linear orderings on
the separate sets $T_{(R,\overline{a})}$.    

\subsection{Some useful lemmas}

We will choose the pair $\mathcal{C}_{1}$, $\mathcal{C}_2$ so that if
$\mathcal{A}$ is $\Delta^{0}_{\alpha}$, then we can make $\mathcal{A}^*$
computable.  We use the lemma below, which was proved in \cite{A-K} (see
also \cite{A-K1}).  We need some definitions.  We write
$(\mathcal{A},\overline{a})\leq_\beta(\mathcal{B},\overline{b})$ if the
$\Pi_\beta$ formulas (of $L_{\omega_1,\omega}$) true of $\overline{a}$ in
$\mathcal{A}$ are true of $\overline{b}$ in $\mathcal{B}$.  If the tuples $\overline{a}$ and $\overline{b}$
are empty, we may write $\mathcal{A}\leq_\beta\mathcal{B}$.  The relations $\leq_\beta$ are called the \emph{standard back-and-forth relations}.  
Let $\alpha$ be a computable ordinal.  A pair
  of structures $\mathcal{C}_{1}$, $\mathcal{C}_{2}$ is said to be
  \emph{$\alpha$-friendly} if the structures are computable and the  
 relations, for $\beta < \alpha$, are
uniformly c.e.  For more about the standard back-and-forth relations, and the notion of $\alpha$-friendliness, see \cite{A-K1}.  
 
  \begin{lem}
\label{lem2.2}
 
  Let $\alpha$ be a computable ordinal, and suppose
$\mathcal{C}_{1}\leq_{\alpha}\mathcal{C}_{2}$,
  where the pair $\{\mathcal{C}_1,\mathcal{C}_2\}$ is $\alpha$-friendly. 
For any $\Pi^0_\alpha$ set $S\subseteq\omega$, there is a uniformly
computable sequence 
  $(\mathcal{A}_{n})_{n\in\omega}$
  such that
  \[\mathcal{A}_{n}\cong\left\{\begin{array}{ll}
  \mathcal{C}_{2} & \mbox{if}\ \ n\in S\\
  \mathcal{C}_{1} & \mbox{otherwise}
  \end{array}\right..\]
 
  \end{lem}
 
  In \cite{A}, there is an analysis of the standard back-and-forth
relations on computable ordinals that yields the next two lemmas.
 
  \begin{lem}
\label{lem2.3}
 
  For a computable ordinal $\beta$, $\omega^{\beta+1}\leq_{2\beta+1}
\omega^{\beta}$ (so $\omega^\beta\leq_{2\beta}\omega^{\beta+1}$), and
there is a computable $\Pi_{2\beta+1}$ sentence that is true in
$\omega^\beta$ but not in $\omega^{\beta+1}$.
 
  \end{lem}
 
  \begin{lem}
\label{lem2.4}
 
  For any computable ordinals $\alpha$, $\beta$, there is an
  $\alpha$-friendly pair of orderings of the types $\omega^{\beta}$,
  $\omega^{\beta+1}$.
 
  \end{lem}       
 
  We need some further definitions and facts from \cite{A-K-M-S},
\cite{C}.  A computable (or
$X$-computable) structure
  $\mathcal{A}$ is said to be \emph{relatively intrinsically
  $\Delta^{0}_{\alpha}$ categorical} if for any copy $\mathcal{B}$ of
  $\mathcal{A}$, there is an isomorphism $f$ from $\mathcal{A}$ onto
  $\mathcal{B}$ such that $f$ is
  $\Delta^{0}_{\alpha}$ relative $\mathcal{B}$ (or $X$ and $\mathcal{B}$).
In \cite{A-K-M-S} and \cite{C}, it is shown that a computable structure 
$\mathcal{A}$ is relatively intrinsically $\Delta^0_\alpha$ categorical iff it
has a \emph{formally $\Sigma^0_\alpha$ Scott family}, where this is a
$\Sigma^0_\alpha$ set of computable $\Sigma_\alpha$ formulas, all with a
fixed tuple of parameters, defining the orbits of all
tuples in $\mathcal{A}$.  If $\mathcal{A}$ has a formally
$\Sigma^0_\alpha$ Scott family with no parameters, then for any copy
$\mathcal{A}$, then we get an isomorphism $f$
which is $\Delta^0_\alpha(\mathcal{B})$ uniformly in
$\mathcal{B}$.  For an $X$-computable structure $\mathcal{A}$, the
result is the same except that the Scott family is a
$\Sigma^0_\alpha(X)$ set of $X$-computable 
$\Sigma_\alpha$ formulas.     

A relation $Q$ is \emph{relatively intrinsically $\Sigma^0_\alpha$} on a
structure $\mathcal{A}$ if in all copies $\mathcal{B}$ of $\mathcal{A}$, 
the image of $Q$ is $\Sigma^0_\alpha$ relative to
$\mathcal{B}$ and $\mathcal{A}$.  In \cite{A-K-M-S} and \cite{C}, it is 
shown that a relation is relatively
$\Sigma^0_\alpha$ on a computable structure
$\mathcal{A}$ iff it is definable by a computable $\Sigma_\alpha$
formula with a finite tuple of parameters in $\mathcal{A}$.     
 
  \begin{lem}
  \label{lem2.5}
 
Suppose $\mathcal{C}_1$ and $\mathcal{C}_2$ are uniformly relatively
$\Delta^0_\beta$ categorical, and each $\mathcal{C}_i$ satisfies a 
computable $\Pi_\beta$ sentence that is not true in the other.  
Let $\mathcal{A}^*$ be as above.     

\begin{enumerate}

\item  If $\mathcal{A}$ is relatively intrinsically
  $\Delta^{0}_{\gamma}$-categorical, then $\mathcal{A}^*$ is relatively
  intrinsically $(\beta+\gamma)$-categorical.  

\item  If $\mathcal{A}^*$ is relatively intrinsically
$\Delta^{0}_{\gamma}$-categorical, then so is
  $\mathcal{A}$.
 
\end{enumerate}

\end{lem}
Proof:  For 1, suppose that $\mathcal{A}$ is relatively intrinsically
  $\Delta^{0}_{\gamma}$-categorical, let\\
  $\mathcal{B}^*\cong\mathcal{A}^*$, and let $\mathcal{B}$ be the copy of
  $\mathcal{A}$ that is naturally defined in $\mathcal{B}^*$, using
computable $\Pi_\beta$ formulas for the basic relations and their
negations.  Then $\mathcal{B}$ is $\Delta^{0}_{\beta}(\mathcal{B}^*)$,
and there is an isomorphism from $\mathcal{A}$ to $\mathcal{B}$ that is
  $\Delta^{0}_{\gamma}(\mathcal{B})$, so it is
$\Delta^0_{\beta+\gamma}(\mathcal{B}^*)$.  We must extend this isomorphism
to the rest of $\mathcal{A}^*$.  For the part $T_{(R,\overline{a})}$,
we locate the corresponding part of $\mathcal{B}^*$, and using
$\Delta^0_\beta(\mathcal{B}^*)$, we can determine whether we are looking
at a pair of copies of $\mathcal{C}_1$ or $\mathcal{C}_2$.  For each
pair of corresponding parts, we get an isomorphism that is
$\Delta^0_\gamma(\mathcal{B}^*)$, uniformly.  Putting the pieces
together, we get an isomorphism from $\mathcal{A}^*$ onto
$\mathcal{B}^*$ that is $\Delta^{0}_{\beta+\gamma}(\mathcal{B}^*)$.
 
For 2, suppose $\mathcal{A}^*$ is relatively intrinsically
  $\Delta^{0}_{\gamma}$-categorical, and let
  $\mathcal{B}\cong\mathcal{A}$.  We may suppose that $\mathcal{B}$ has
universe equal to the set of even numbers, and we can construct
  $\mathcal{B}^*\cong\mathcal{A}^*$ over the given $\mathcal{B}$, such
that $\mathcal{B}^*\leq_T\mathcal{B}$.
  For an isomorphism from $\mathcal{A}^*$ onto $\mathcal{B}^*$ that is
  $\Delta^{0}_{\gamma}(\mathcal{B}^*)$, the restriction is an isomorphism
from $\mathcal{A}$ onto $\mathcal{B}$ that is
$\Delta^0_\gamma(\mathcal{B})$.  Therefore,
  $\mathcal{A}$ is relatively intrinsically
  $\Delta^{0}_{\gamma}$-categorical.

  \begin{lem}
  \label{lem2.6}
 
Form $\mathcal{A}^*$ as above, using $\mathcal{C}_1$ and
$\mathcal{C}_2$ which are relatively $\Delta^0_\beta$ categorical,
uniformly, where each satisfies a computable $\Pi_\beta$ sentence
that is not true in the other.  Let $Q$ be a relation on $\mathcal{A}$. 

\begin{enumerate}

\item  If $Q$ is relatively intrinsically $\Sigma^{0}_{\gamma}$
on $\mathcal{A}$, then it is relatively intrinsically
$\Sigma^{0}_{\beta+\gamma}$ on $\mathcal{A}^*$.  

\item

If $Q$ is relatively intrinsically
  $\Sigma^{0}_{\gamma}$ on $\mathcal{A}^*$, then it is relatively
  intrinsically $\Sigma^{0}_{\gamma}$ on $\mathcal{A}$.

\end{enumerate}
 
  \end{lem}
 
  The proof of Lemma \ref{lem2.6} is like that for Lemma \ref{lem2.5}.
 
  \bigskip

Suppose $\mathcal{A}$ is $\Delta^0_\alpha$.  We may assume that
the universe is computable.  Take $\beta$ such that $\alpha\leq
2\beta+1$, let $\mathcal{C}_1$, $\mathcal{C}_2$ be orderings of types
$\omega^{\beta+1}$, $\omega^{\beta}$, respectively, and form
$\mathcal{A}^*$ as described.  By Lemmas \ref{lem2.3},
\ref{lem2.4}, and \ref{lem2.5}, we may take $\mathcal{A}^*$ to be
computable.  The orderings $\mathcal{C}_1$, $\mathcal{C}_2$ are both uniformly relatively $\Delta^0_{2\beta+2}$ categorical.  Since there is a computable $\Pi_{2\beta+1}$ sentence true in $\omega^\beta$ but not in $\omega^{\beta+1}$, each $\mathcal{C}_i$ 
satisfies a computable $\Pi_{2\beta+2}$ sentence that is not true in the
other.  There is a $\Delta^0_{2\beta+2}$ copy of
$\mathcal{A}$ definable in $\mathcal{A}^*$---with an atomic formula
defining the universe $A$, and computable $\Pi_{2\beta+2}$ formulas
defining the basic relations of $\mathcal{A}$ and their negations.

\subsection{Proof of Theorem \ref{thm2.1}}

We are ready to show that the Scott ranks of $\mathcal{A}$ and
$\mathcal{A}^*$ are related as in Theorem \ref{thm2.1}.  If the orbit of a
tuple in one of the structures is
  intrinsically hyperarithmetical, then by a result of
  Soskov \cite{S}, it is definable by a computable infinitary formula
  without parameters (see also \cite{G-H-K-S}).  Suppose one of the
structures has computable rank.
  Then that structure is relatively intrinsically
  hyperarithmetically categorical. By Lemma \ref{lem2.6}, the second
structure is also relatively intrinsically
  hyperarithmetically categorical, so for some computable $\gamma$, there
  is a formally $\Sigma^{0}_{\gamma}$ Scott family, possibly with
parameters.  The orbit of the parameters is itself relatively
intrinsically
  hyperarithmetical, so it is defined by a computable infinitary formula
  without parameters.  Therefore, we can remove the parameters in the
  Scott family.  From this, it follows that the second structure has 
computable Scott rank.
 
  Next, suppose $\mathcal{A}$ has rank $\omega_{1}^{CK}$.  The
  orbits of all tuples in $\mathcal{A}$ are definable by computable
  infinitary formulas.  Consider a tuple in $\mathcal{A}^*$.  The part in
  $\mathcal{A}$ is described by a computable infinitary formula, and the
  part outside is as well.  So, the orbits are all defined by
  computable infinitary formulas in $\mathcal{A}^*$.  Suppose
  $\mathcal{A}^*$ has rank $\omega_{1}^{CK}$.  The orbits of all tuples in
  $\mathcal{A}^*$ are defined by computable infinitary formulas.  In
  $\mathcal{A}$, the orbits are all relatively intrinsically
  hyperarithmetical, so they are defined by computable infinitary 
  formulas,
  without parameters.
 
From what we have said, it follows that if one
  of the structures has rank $\omega_{1}^{CK}+1$, then
  the other does as well.  
This completes the proof Theorem \ref{thm2.1}.  Combining this
with Makkai's Theorem (Theorem
\ref{thm1.3}), we obtain a computable structure of rank $\omega_1^{CK}$, as in Theorem  \ref{thm1.4}.     

\section{From trees to structures}  

We have given one proof of Theorem \ref{thm1.4}, taking Theorem
\ref{thm1.3} as given.  We shall re-work the proof of Theorem
\ref{thm1.3} to give another, self-contained, proof of Theorem
\ref{thm1.4}.  In this section, we reduce the proof to the construction 
of a computable tree with special properties.  The construction is essentially 
due to Makkai.  Morozov \cite{Mo} used it to produce computable 
structures which share with the Harrison ordering the feature that there 
are non-trivial automorphisms, but not hyperarithmetical ones.  We follow
Morozov in the details.         

Let $T$ be a subtree of $\omega^{<\omega}$, and let $T_n$ 
be the set of elements of $T$ of length $n$.  Let $G_n$ 
be the set of finite subsets of $T_n$.  For each $n$, $G_n$ is an Abelian group
under the operation of symmetric difference, denoted by $\triangle$.  
The identity is $\emptyset$, and each element is its own inverse.  In fact, 
$G_n$ looks like a vector space over the $2$-element 
field, with a basis consisting of $\{t\}$, for $t\in T_n$.  For
distinct $n > 0$, we replace $\emptyset$ by distinct elements 
$id_n$---$G_0$ consists of the
single element $\emptyset = id_0$.    

Let $G = \cup_n G_n$.  The tree structure on $T$ induces a tree 
structure on $G$, which is most easily described in terms of a
 \emph{predecessor} function, which we will denote by $p$.  If $a\in G_{n+1}$ 
represents a non-empty subset of $T_{n+1}$, say
$a = \{t_1,\ldots,t_k\}$, where $t_i$ has predecessor $t_i'$ (in the sense of $T$), 
then $p(a)$ is the sum (in $G_n$) of the elements $\{t_i'\}$.  We let $p(id_{n+1}) = id_n$.   If $t\in T_n$ has two successors $t',t''\in T_{n+1}$, then 
$p(\{t',t''\}) = id_n$.  For each $n$, $p$ acts as a homomorphism from $G_{n+1}$ into $G_n$.  Using the group structures, together with the tree structure, we define a 
family of unary operations $f_a$, for $a\in G$.  If $a\in G_n$, then for $b\in G_m$, let 
$k = min\{m,n\}$, and let $a^*$ and $b^*$ be the result of iterating $p(x)$ until the 
value is in $G_k$ (i.e., $a^* = p^{(n-k)}(a)$ and $b^* = p^{(m-k)}(b)$).  
Then $f_a(b) = a^*\triangle b^*$---this is the same as $f_b(a)$.  Finally, we let 
$\mathcal{A}(T)$ be the structure $(G,(f_a)_{a\in G})$.

We have $a\in G_n$ iff $f_{id_n}(a) = a$ and for all $m < n$, $f_{id_m}(a)\not= a$.  It follows that $G_n$ is preserved under automorphisms of $\mathcal{A}(T)$.  Similarly, $p$ is preserved under automorphisms, since for $a\in G_{n+1}$ and $b\in G_n$, we have $p(a) = b$ iff $f_{id_n}(a) = b$.  For any $a\in G_n$, $a = f_a(id_n)$.  It follows that any automorphism of $\mathcal{A}(T)$ is determined by the images of the elements $id_n$.     

\begin{thm} [Morozov]

For any computable tree $T\subseteq\omega^{<\omega}$, $\mathcal{A}(T)$ is a computable structure such that $T$ has a path iff $\mathcal{A}(T)$ has a non-trivial automorphism.  Moreover, if $T$ has paths, then the Turing degrees of of the paths are the same as the degrees of non-trivial automorphisms of $\mathcal{A}(T)$.  

\end{thm} 
Proof:  Clearly, $\mathcal{A}(T)$ is computable in $T$.  Suppose $g$ is a non-trivial automorphism, and let $g_n = g(id_n)$.  Since $g$ is non-trivial, there is some $n$ such that $g_n\not= id_n$.  Take $t_n\in g_n$.  For $m < n$, let
$t_m$ be the level $m$ predecessor of $t_n$ in $T$.  By induction, for $m > n$, choose $t_m$ to be a successor of $t_{m-1}$ in $g_m$.  These elements form a path through
$T$.  Conversely, if $(t_n)_{n\in\omega}$ is a path, we obtain a non-trivial automorphism $g$ by mapping $id_n$ to $\{t_n\}$, and then letting $g(a) = f_a(g(id_n))$, for $a\in G_n$.

\bigskip            

To obtain a computable structure of rank $\omega_1^{CK}$, we will produce a computable tree $T$ with special features, and form the structure $\mathcal{A}(T)$.  We need to develop an understanding of Scott rank in these structures.   

\subsection{Tree rank}

We assign ranks to nodes of a tree $T$, and to $T$ itself, as follows.  
\begin{enumerate}

\item  For $\sigma\in T$, $rk(\sigma) = 0$ if
$\sigma$ has no successors.  

\item  For $\alpha > 0$, $rk(\sigma) = \alpha$ if
$\alpha$ is the first ordinal greater than $rk(\tau)$ for all successors
$\tau$ of $\sigma$.

\end{enumerate}
Note that $rk(\sigma)\geq\alpha+1$ iff $\sigma$ has a successor of rank $\geq\alpha$, and for limit $\alpha$, $rk(\sigma)\geq\alpha$ iff $rk(\sigma)\geq\beta$ for all $\beta <\alpha$.  We write $rk(\sigma) = \infty$ if $\sigma$ does not
have ordinal rank, where we consider that $\infty \geq \alpha$, for all
ordinals $\alpha$.  We assign to the tree $T$ the rank of the top node; i.e., $rk(T) = rk(\emptyset)$. 

\begin{lem}
\label{lem3.2}

For all $\sigma\in T$, $rk(\sigma) = \infty$ iff $\sigma$ lies on a path.

\end{lem}

We connect Scott rank in $\mathcal{A}(T) = (G,(f_a)_{a\in G})$ to tree rank in $T$.  
This will require several steps.  First, for $a\in G$, we calculate the rank of $a$ in the derived tree structure on $G$ in terms of the ranks of elements of $a$ in $T$.

\begin{lem}
\label{lem3.3}

Let $a\in G_n$, where $a\not= id_n$.  Then $rk(a)$ is the minimum of $rk(t)$, for $t\in a$.

\end{lem}
Proof:  Say $a = \{t_1,\ldots,t_k\}$.  
We show, by induction on $\alpha$, that $rk(a) \geq \alpha$ iff for all $t\in a$, $rk(t)\geq \alpha$.  The statement is clearly true for $\alpha = 0$, and the clause for limit $\alpha$ is trivial.  We consider $\alpha = \beta + 1$.  First, suppose 
$rk(a) \geq \beta + 1$.  Then $a$ has a successor $a^*$ of rank $\geq\beta$.  By H.I., all 
$t\in a^*$ have rank $\geq\beta$.  Now, 
$a^*= \{t_1^*,\ldots,t_2^*\}\cup k$, 
where $t_i^*$ is a successor of $t_i$ of least rank in 
$a^*$, and $k$ is a successor of $id_n$
(if $t_i$ has three successors in $a^*$, then 
$t_i^*$ is one of least rank, and the other two are in $k$).  
Since $t_i^*$ has rank $\geq\beta$, $t_i$ has rank $\geq\beta+1$.  For the other direction, suppose that all $t_i$ have rank $\geq\beta +1$.  For each $i$, let $t_i^*$ be a successor of $t_i$ of rank $\geq\beta$.  Then $a^* = \{t_i^*,\ldots,t_k^*\}$ is a successor of $a$, and by H.I., $rk(a^*)\geq\beta$.  Therefore, $rk(a)\geq\beta + 1$. 

\bigskip

From any single element $a\in G_n$, we can define any $b\in \cup_{m\leq n}G_m$, by the atomic formula $x = f_b(f_a(a))$.  Using this, we see that if $a,a^*\in G_n$, then the function taking $a$ to $a^*$ extends in a unique way to an automorphism of the substructure $\cup_{m\leq n}G_m$,           

\begin{lem}
\label{lem3.4}

Let $\overline{a}$ and $\overline{b}$ be tuples such that
$\overline{a}\equiv^0\overline{b}$.  Suppose $n$ is greatest such that
some $a_i$ is in $G_n$, and let $g = f_{a_i}(b_i)$.  Then for all $\beta$, 
the following are equivalent:

\begin{enumerate}

\item  $\overline{a} \equiv^\beta \overline{b}$,

\item  $\overline{a},id_n \equiv^\beta \overline{b},g$,

\item  $id_n\equiv^\beta g$.

\end{enumerate}

\end{lem}
Proof:  We proceed by induction.  For $\beta = 0$, Statements 1, 2, and 3 are all true.  The clause for limit $\beta$ is trivial.  We consider $\beta = \gamma + 1$.  Supposing 1, we prove 2.  For any $c$, there exists $d$, and for any $d$, there exists $c$, such that $\overline{a},c\equiv^\gamma\overline{b},d$.  We suppose that $c$ and $d$ are in $G_m$ for $m > n$---this is the interesting case.  Let 
$g' = f_c(d)$.  By H.I., $\overline{a},c,id_n,id_m\equiv^\gamma\overline{b},d,g,g'$.  Therefore, $\overline{a},id_n,c\equiv^\gamma\overline{b},g,c$.  This gives 2.
Statement 2 clearly implies Statement 3.  Supposing 3, we prove 1.  For any $c$, there exists $d$, and for any $d$, there exists $c$ such that $id_n,c\equiv^\gamma a,d$.  As above, we suppose that $c,d\in G_m$ for $m > n$ and let $g' = f_c(d)$.  By H.I., we get $\overline{a},c\equiv^\gamma\overline{b},d$.  This gives 1.  

\bigskip

From Lemma \ref{lem3.4}, we get the following.

\begin{lem}
\label{lem3.5}

For any tuple $\overline{a}$, and $n$ is greatest such
that some $a_i\in G_n$, then $SR(\overline{a}) = SR(id_n)$.   

\end{lem}

The Scott rank of $id_n$ is calculated as follows.  

\begin{lem}
\label{lem3.6} 

For $a\in G_n$, $a \equiv^{\beta} id_n$ iff $rk(a)\geq\omega\cdot\beta$.
  
\end{lem}
Proof:  We proceed by induction on $\beta$.  

\bigskip
\noindent
\textbf{Case 1}:  $\beta = 0$.

\bigskip

We have a unique automorphism of the substructure $G_{\leq n} = \cup_{m\leq n}G_m$ taking $id_n$ to $a$.  It follows that for any atomic formula $\varphi(x)$, 
$\varphi(id_n)$ holds iff $\varphi(a)$ holds.       

\bigskip
\noindent
\textbf{Case 2}:  $\beta = \gamma + 1$.

\bigskip

First, we suppose that $a \equiv^{\gamma+1} id_n$, and show that $rk(a)\geq
\omega(\gamma+1)$.  For each $k\in\omega$, there exists 
$c\in G_{n+k}$ such that $a,c\equiv^\gamma id_n,id_{n+k}$.  Then $c\equiv^\gamma id_{n+k}$, and by H.I., $rk(c)\geq\omega\gamma$.  Since $f_{id_n}(id_{n+k}) = id_n$, we must have 
$f_{id_n}(c) = a$, so $c$ is a $k^{th}$ successor of $a$.  Therefore, 
$rk(a)\geq\omega\gamma+k$.  Since this is true for all $k$, 
$rk(a)\geq\omega(\gamma+1)$.

Now, we suppose that $rk(a)\geq\omega(\gamma+1)$ and show that
$a\equiv^{\gamma+1} id_n$.  Suppose $b\in G_m$, where $m = n+k$ (this is the interesting case).  Take $a^*\in G_m$ a $k^{th}$ successor of $a$ such that $rk(a^*)\geq\omega\gamma$.  By H.I., $a^*\equiv^\gamma id_m$.  Letting $c = f_b(a^*)$, we have $a^*,a,c\equiv^0 id_m,id_n,b$.  Similarly, letting $c' = f_b(f_{a^*}(id_m)$, we have $a^*,a,b\equiv^0 id_m,id_n,c'$.  By Lemma \ref{lem3.4}, 
we can replace $\equiv^0$ by $\equiv^\gamma$ in both cases.  Then we get 
$a,c\equiv^\gamma id_n,b$ and $a,b\equiv^\gamma id_n,c'$.  Therefore, $a\equiv^{\gamma+1} id_n$.   

\bigskip
\noindent
\textbf{Case 3}:  $\beta$ is a limit ordinal.  

\bigskip

In this case, the statement holds by induction.
           
\bigskip

Having completed the proof of Lemma \ref{lem3.6}, we come to the main result of the section.       

\begin{thm}
\label{thm3.7}  

Suppose $T$ is a computable tree such that for each $n$, there is a
computable bound on the ordinal ranks of elements of
$T_n$, but there is no computable bound on the set of ordinal ranks of
elements of $T$.  Then $\mathcal{A}(T)$ is computable, and
$SR(\mathcal{A}(T)) = \omega_1^{CK}$.

\end{thm}
Proof:  It is enough to show that for each $n$, $SR(id_n) <
\omega_1^{CK}$, but there is no computable upper bound on these ranks.  We
have $SR(id_n) > \beta$ iff there exists $a$ not in the orbit of $id_n$
such that $rk(a)\geq\omega\beta$.  For each $n$, there is a
computable bound on the ordinal ranks of elements of $G_n$, say
$\omega\beta_n$ is greater, but still computable.  If $a\in G_n$ and
$rk(a)\geq\omega\beta_n$, then $a$ is in the orbit of $id_n$, so
$SR(id_n) \leq \beta_n$.  For each computable $\beta$, there exists
$n$ and $a\in G_n$ such that $rk(a)$ is an ordinal $\geq\omega\beta$. 
Then $a\equiv^\beta id_n$, but $a$ is not in the orbit of $id_n$. 
This completes the proof of Theorem \ref{thm3.7}. 

\bigskip

So, to produce a computable structure of rank $\omega_1^{CK}$, 
it is enough to produce an appropriate computable tree, with a computable bound on the ordinal ranks at each level, but no computable bound over-all.  

\section{Thin trees}

We define a special kind of tree satisfying the
conditions of Theorem \ref{thm3.7}.

\bigskip
\noindent
\textbf{Definition}:  Let $T$ be a subtree of $\omega^{<\omega}$, and let
$T_n$ be the set of nodes at level $n$.  We say that $T$ is \emph{thin}
if for each $n\geq 1$, the order type of the set of ordinal ranks of nodes in
$T_n$ is $\leq\omega\cdot n$.

\begin{lem}
\label{lem4.2}
  
If $T$ is a computable thin tree, then for each
$n$, there is a computable bound on the ordinal ranks of elements of
$T_n$.

\end{lem}

Sacks suggested a way to produce computable thin trees of arbitrarily
high computable rank.

\begin{lem}
\label{lem4.3}  

For each computable ordinal $\alpha$, there is a computable
thin tree $T$, with no paths, such that $rk(T)\geq\alpha$.

\end{lem}  
Proof:  Let $a$ be a notation for $\alpha$ (in Kleene's $\mathcal{O}$).  We define a preliminary computable tree $T_a$, with labels, as follows.  We give $\emptyset$ the label $a$.  If $\sigma$ has label $1$ (representing $0$), then $\sigma$ has no successors.  If $\sigma$ has label $b = 2^c$, then $\sigma$ has a single successor, with label $c$.  If $\sigma$ has label $3\cdot 5^e$, then $\sigma$ has an infinite family of successors, with labels $b_n$, where $\varphi_e(n) = b_n$.  If $\sigma\in T$ has label $b$, where $|b| = \beta$, then $\beta$ is the rank of $\sigma$ in $T_a$.  Then $T_a$ has rank $\alpha$.  We replace $T_a$ by a thin tree by slowing down the expansion of the limit notations.  

We define the tree $T$, together with a computable function $t$ from $T_a$ $1-1$ into $T$.  For $\sigma$ at level $0$ or $1$ in $T_a$, $t(\sigma) = \sigma$.  If $\sigma\in T$ has label of form $2^b$, and $\sigma'$ is the unique successor, with label $b$, then $t(\sigma)$ has a unique successor in $T$, and this is $t(\sigma')$.  If $\sigma\in T_a$ has a label of form $3\cdot 5^e$, where $\varphi_e(n) = b_n$, then we regard $\sigma$ as ``needing expansion''.  We may delay expansion, in which case, $t(\sigma)$, will have a single successor, not in $ran(t)$.  We may delay several times, giving $t(\sigma)$ a finite chain of successors.  Say $\sigma$ has successors $\sigma_n$ with labels $b_n$.  When we expand $\sigma$, we define $t(\sigma_n)$ to be a successor of the element $\sigma'$ of length $s$ in the chain below $t(\sigma)$.        

At stage $s$ in the construction, we have defined $t$ mapping a finite subtree $T_a^s$ of $T_a$ into a finite subtree $T^s$ of $T$, in which all terminal elements have length $s$ or  are $t$-images of elements with label $1 = |0|$.  The subtree $T_a^s$ may have some nodes $\sigma$ which need expansion.  If so, then at stage $s+1$, we choose the first such $\sigma$, not previously acted on, and we begin the expansion process, letting $t$ map the first successor of $\sigma$ to a new successor of the level $s$ extension of $t(\sigma)$.  For other nodes $\sigma\in T_a^s$ such that $t(\sigma)$ has length $s$ in $T$, we add a single successor, and we extend $t$, if appropriate.  Finally, for each of the expansions begun earlier, we continue the expansion process.  If we are expanding $\sigma$, then we extend $t$, mapping the next successor $\sigma_n$ in $T_a - T_a^s$ to a new successor of $t(\sigma)$.  In addition, we extend $t$ further, and we add further nodes to $T$, so that the terminal nodes will all be at level $s+1$ or be $t$-images of elements with label $1$.  The tree $T$ obtained in this way is clearly computable.       

It is easy to see that if $\sigma\in T$ has rank $\beta$, then $t(\sigma)$ has rank $\geq\beta$, so $T$ has rank $\geq\alpha$.  We can show, by induction on $\beta$, that if $\sigma\in T$ has rank $\beta$, then $t(\sigma)$ has rank $< \beta + \omega$.  
We must show that $T$ is thin.  We show by induction on $m\geq 1$ that the set of ordinal ranks of elements at level $m$ is at most $\omega\cdot m$.      
At level $1$, in $T_a$, there is either just one element or an infinite sequence $\sigma_n$  with increasing ranks $\beta_n$.  In the second case, note that for each $n$, $t(\sigma_n) = \sigma_n$, and the rank of $\sigma_n$ in $T$ is $\beta_n+k_n$, for some $k_n <\omega$.  The sequence may not be strictly increasing, but for each $n$, there are only finitely many $r$ such that $\beta_r+k_r < \beta_n+k_n$, so the set of ranks still has order type $\omega$.  

Supposing that the set of ranks at level $m$ in $T$ has order type at most $\omega\cdot m$, we show that the set of ranks at level $m+1$ has order type at most $\omega\cdot (m+1)$.  There is at most one node at level $m$ in $T$ with infinitely many successors at level $n+1$.  All other nodes at level $m$ have at most one successor.  Say $\sigma$ is the element of $T_a$ chosen for expansion.  Suppose $\sigma$ has limit rank $\beta$, with successors $\sigma_n$ of rank $\beta_n$.  Let $\tau$ be the extension of $t(\sigma)$ at level $m$ in $T$.  The successors of $\tau$ are the elements $\tau_n = t(\sigma_n)$, where $\beta_n\leq rk(\tau_n) < rk(\beta)$.  Therefore, $\tau$ has rank $\beta$.  While $rk(\tau_n)$ need not be strictly increasing with $n$, since $rk(\tau_n)$ differs only finitely from $\beta_n$, for each $n$, there are only finitely many $r$ such that $rk(\tau_r) < rk(\tau_n)$.  For the elements at level $m$ in $T$ aside from $\tau$, either the rank is $0$, and there is no successor, or the rank is $\gamma+1$, for some $\gamma$, and the unique successor has rank $\gamma$.  So, at level $m+1$ in $T$, the set of ranks is obtained by reducing the successor ranks that occur at level $m$ by $1$, and then adding an increasing sequence of ranks.  The order type for level $m+1$ is the result of taking $m$ copies of $\omega$, or an initial segment, and adding an $\omega$ sequence.  The result has order type at most $\omega\cdot (m+1)$.         

\begin{lem}
\label{lem4.4}  

There is a computable thin tree $T$ such that $rk(T) = \infty$, but $T$ has no hyperarithmetical path.

\end{lem}
Proof:  We apply Barwise-Kreisel Compactness, using the previous lemma to
satisfy the hypotheses.  We have a $\Pi^1_1$ set $\Gamma$ of
computable infinitary sentences describing a model of $KP$ in which the
ordinals have an initial segment of type $\alpha$ for each computable
ordinal $\alpha$.  We use a unary relation symbol $T$, and include sentences saying that
it is a computable tree that does not have computable ordinal rank, but
has no hyperarithmetical path.  We add a binary relation symbol $F$, and include axioms
saying that $F$ maps elements of $T$ to ordinals such that for $t\in T$
and $\alpha$ a computable ordinal, $F(t) = \alpha$ iff
$rk(t) = \alpha$.  Also, for each $n$, we add a binary relation symbol
$C_n$, and we include axioms saying that $C_n$ is an order-preserving
function from $ran(F|T_n)$ onto an initial segment of $\omega\times n$.   

Any $\Delta^1_1$ subset of $\Gamma$ is satisfied by taking the 
least admissible set as the model of $KP$, letting $T$ be as in Lemma
\ref{lem4.3}, of sufficiently high computable rank, with no path, 
taking $F$ to be the rank function, and letting $C_n$ be the
function collapsing the ranks of elements of $T_n$ to an initial segment
of $\omega\cdot n$.  Therefore, there is a model of the whole of
$\Gamma$.  In this model, the ordinals will not be well-founded.  The
function $F$ will map tree elements that have ordinal rank to their
ranks, while mapping unranked elements to non-standard ordinals.  The
function $C_n$ will map the ordinal ranks of elements of $T_n$ that have
ordinal rank to an initial segment of $\omega\cdot n$.           

\bigskip

The fact that the tree in Lemma \ref{lem4.4} does not have computable
ordinal rank and has no hyperarithmetical path implies that there is no
computable bound on the ordinal ranks of nodes.  So, we get a tree $T$ as
in Theorem \ref{thm3.7}, and then $\mathcal{A}(T)$ is the structure we
need for Theorem \ref{thm1.4}.

\section{Categoricity}

The Ryll-Nardjewski Theorem says that for an elementary first order theory
$T$ which is $\aleph_0$ categorical, there are only finitely many
complete types in any fixed tuple of variables.  Each type is principal,
so there is a finite set of generating formulas.  It is natural to
wonder whether this generalizes.  Suppose $T$ is a computable infinitary
theory which is $\aleph_0$ categorical.  The types realized in models of
$T$ must be principal; i.e., for each type, there is a computable
infinitary formula that generates it.  However, for a given tuple of
variables, there need not be a nice set of formulas generating the
types.  For the structures $\mathcal{A}(T)$ that we have been considering,
the Ryll-Nardjewski Theorem fails badly.      

\begin{prop}

If $T$ is a computable thin tree which is unranked but has no 
hyperarithmetical path, then the computable
infinitary theory of $\mathcal{A}(T)$ is $\aleph_0$ categorical.  If
$\Phi$ is a set of computable infinitary formulas defining the orbits
of elements of $\mathcal{A}(T)$, then $\Phi$ cannot be
hyperarithmetical, or even $\Pi^1_1$.  

\end{prop}
Proof:  For each $n$, $id_n$ has computable Scott rank, and the orbit is
defined by a computable infinitary formula $\varphi_n(u)$.  For any $a\in G_n$, the
orbit of $a$ is defined by $\exists u\,\varphi_n(u)\ \&\ x = f_a(u)$.
There is no hyperarithmetical set of computable infinitary formulas
defining the orbits of all elements, for then there would be a computable
bound on the ranks.  It follows that there can be no
$\Pi^1_1$ set $\Phi$ of computable infinitary formulas defining the orbits
of all elements, for then we could apply Barwise-Kreisel Compactness to
produce a model of the computable infinitary theory of
$\mathcal{A}(T)$ with an element $a$ satisfying
$\neg{\varphi(x)}$, for all $\varphi(x)\in\Phi$.   

\bigskip

\section{Problems}

We close with some open problems. 

\bigskip

The Harrison ordering $\mathcal{H}$, and the Harrison Boolean algebra the Harrison groups, have a strong computable approximation property.     

\begin{prob}

Is there a structure $\mathcal{A}$ of rank $\omega_1^{CK}$ such that for any $\Sigma^1_1$ set $S$, there is a uniformly computable sequence $(\mathcal{A}_n)_{n\in\omega}$ such that if $n\in S$, then $\mathcal{A}\cong\mathcal{A}$ and if $n\notin S$, then $\mathcal{A}_n$ has computable rank?       

\end{prob}

\bigskip

Sacks and Young have further examples of structures of rank
  $\omega_{1}^{CK}$ with properties different from Makkai's example.
  Their examples are not hyperarithmetical but are ``hyperarithmetically
  saturated'', which means that the ordinals in the least
  admissible set over the structure are just the computable ordinals.
 
  \begin{prob}

Can we convert hyperarithmetically saturated structures of rank $\omega_{1}^{CK}$ into computable structures of this rank?

\end{prob}


\begin{thebibliography}{99}
 
\bibitem{A}  Ash, C.\ J., ``Categoricity in hyperarithmetical
degrees'', \emph{Annals of Pure and Applied Logic}, vol.\ 34(1987), 
pp.\ 1--14.
 
\bibitem{A-K}  Ash, C. J., and J. F. Knight, ``Pairs of recursive
structures'', \emph{Annals of Pure and Applied Logic}, vol. 46(1990), 
pp.\ 211-234.

\bibitem{A-K1}  Ash, C.\ J., and J.\ F.\ Knight, \emph{Computable
structures and the hyperarithmetical hierarchy}, Elsevier, 2000.

\bibitem{A-K-M-S} Ash, C.\ J., J.\ F.\ Knight, M.\ Mannasse, and T.\
Slaman, ``Generic copies of countable structures'', \emph{Annals of
Pure and Applied Logic}, vol.\ 42(1989),  pp.\ 195--205.

\bibitem{C}  Chisholm, J., ``Effective model theory versus recursive
model theory'',\linebreak \emph{J.\ Symb. Logic}, vol.\ 55(1990), pp.\
1168--1191. 

\bibitem{G-H-K-S}  Goncharov, S.\ S., V.\ S.\ Harizanov, J.\ F.\ Knight,
and R.\ I.\ Shore, ``Intrinsically $\Pi^1_1$ relations and paths through
$\mathcal{O}$'', submitted to \emph{J.\ Symb. Logic}. 

\bibitem{G-H-K-M-M-S}  Goncharov, S.\ S., V.\ S.\ Harizanov, J.\ F.\
Knight, C.\ McCoy, R.\ G.\ Miller, and R.\ Solomon, ``Enumerations in
computable structure theory'', submitted to \emph{Annals of Pure and
Applied Logic}.

\bibitem{H} Harrison, J., ``Recursive pseudo
well-orderings'',\emph{Transactions of the Amer. Math. Soc.}, vol.\
131(1968), pp.\ 526--543,
 
\bibitem{M}  Makkai, M., ``An example concerning Scott heights'',
  \emph{J.\ Symb.\ Logic}, vol.\ 46(1981), pp.\ 301--318.

\bibitem{Mo}  Morozov, A.\ S., ``Groups of computable automorphisms'',
in \emph{Handbook of Recursive Mathematics}, ed.\ by Yu.\ L.\ Ershov,
et.\ al., 1998, Elsevier, pp.\ 311--345.

\bibitem{Sa} Sacks, G.\ E., \emph{Higher Recursion Theory},
   Berlin:Springer-Verlag, pp.\ 8--18.

\bibitem{Sc}  Scott, D., Scott, D., ``Logic with denumerably long 
formulas and finite strings of quantifiers'', in \emph{The Theory of
Models}, ed.\ by J.\ Addison, L.\ Henkin, and A.\ Tarski, North-Holland,
1965, pp.\ 329--341.
 
\bibitem{S}  Soskov, I.\ N., ``Intrinsically hyperarithmetical sets'',
   \emph{Math. Logic Quarterly}, vol.\ 42(1996), pp.\ 469--480.

\end{thebibliography}
  \end{document}